\documentclass[11pt,leqno]{amsart}
\usepackage{amssymb,verbatim,enumerate,ifthen,times}
\usepackage[mathscr]{eucal}
\usepackage[utf8]{inputenc}
\usepackage[T1]{fontenc}
\usepackage{tabu}
\usepackage{dsfont}
\usepackage{accents}

\def\N{\mathbb{N}}

\def\Q{\mathbb{Q}}
\def\R{\mathbb{R}}

\def\A{\mathscr{A}}
\def\M{\mathscr{M}}

\def\id{\mathop{\mbox{\rm id}}\nolimits}

\def\Ref(#1|#2){(#1\hspace{.3mm}|\hspace{.3mm}#2)}
\newtheorem{theorem}{Theorem}
\newtheorem*{theorem*}{Theorem}
\long\def\Thm#1#2{\ifthenelse{\equal{#1}{*}}{\begin{theorem*}#2\end{theorem*}}
             {\begin{theorem}\label{T#1}#2\end{theorem}}}
\newtheorem{Atheorem}{Theorem}

\def\thm#1{Theorem~\ref{T#1}}

\newtheorem{proposition}[theorem]{Proposition}
\newtheorem*{proposition*}{Proposition}
\long\def\Prp#1#2{\ifthenelse{\equal{#1}{*}}{\begin{proposition*}#2\end{proposition*}}
             {\begin{proposition}\label{P#1}#2\end{proposition}}}
\def\prp#1{Proposition~\ref{P#1}}

\newtheorem{corollary}[theorem]{Corollary}
\newtheorem*{corollary*}{Corollary}
\long\def\Cor#1#2{\ifthenelse{\equal{#1}{*}}{\begin{corollary*}#2\end{corollary*}}
             {\begin{corollary}\label{C#1}#2\end{corollary}}}

\newtheorem{lemma}[theorem]{Lemma}
\newtheorem*{lemma*}{Lemma}
\long\def\Lem#1#2{\ifthenelse{\equal{#1}{*}}{\begin{lemma*}#2\end{lemma*}}
             {\begin{lemma}\label{L#1}#2\end{lemma}}}
\def\lem#1{Lemma~\ref{L#1}}

\theoremstyle{definition}
\newtheorem{definition}[theorem]{Definition}
\newtheorem*{definition*}{Definition}
\long\def\Defn#1#2{\ifthenelse{\equal{#1}{*}}{\begin{definition*}\rm #2\end{definition*}}
             {\begin{definition}\label{D#1}\rm #2\end{definition}}}

\newtheorem{remark}[theorem]{Remark}
\newtheorem*{remark*}{Remark}
\long\def\Rem#1#2{\ifthenelse{\equal{#1}{*}}{\begin{remark*}\rm #2\end{remark*}}
             {\begin{remark}\label{R#1}\rm #2\end{remark}}}
\def\rem#1{Remark~\ref{R#1}}

\newtheorem{example}{Example}
\newtheorem*{example*}{Example}
\long\def\Exa#1#2{\ifthenelse{\equal{#1}{*}}{\begin{example*}\rm #2\end{example*}}
             {\begin{example}\label{Ex#1}\rm #2\end{example}}}

\def\eq#1{{\rm(\ref{E#1})}}
\def\Eq#1#2{\ifthenelse{\equal{#1}{*}}
  {\begin{equation*}\begin{aligned}[]#2\end{aligned}\end{equation*}}
  {\begin{equation}\begin{aligned}\label{E#1}#2\end{aligned}\end{equation}}}

\begin{document}
\vspace{5mm}

\date{\today}

\title{On the balancing property of Matkowski means}
\author[T. Kiss]{Tibor Kiss}
\address{
\textit{Address of the author:}
Institute of Mathematics, University of Debrecen, H-4002 Debrecen, Egyetem tér 1, Hungary}
\email{kiss.tibor@science.unideb.hu}
\thanks{
The research of the author was supported in part by the NKFIH Grant K-134191 and in part by the project no. 2019-2.1.11-TÉT-2019-00049, which has been implemented with the support provided from the National Research, Development and Innovation Fund of Hungary, financed under the TÉT funding scheme.}

\begin{abstract}
Let $I\subseteq\R$ be a nonempty open subinterval. We say that a two-variable mean $M:I\times I\to\R$ enjoys the \emph{balancing property} if, for all $x,y\in I$, the equality
\Eq{GA}{
M\big(M(x,M(x,y)),M(M(x,y),y)\big)=M(x,y)
}
holds.

The above equation has been investigated by several authors. The first remarkable step was made by Georg Aumann in 1935. Assuming, among other things, that $M$ is \emph{analytic}, he solved \eq{GA} and obtained quasi-arithmetic means as solutions. Then, two years later, he proved that \eq{GA} characterizes \emph{regular} quasi-arithmetic means among Cauchy means, where, the differentiability assumption appears naturally. In 2015, Lucio R. Berrone, investigating a more general equation, having symmetry and strict monotonicity, proved that the general solutions are quasi-arithmetic means, provided that the means in question are \emph{continuously differentiable}.

The aim of this paper is to solve \eq{GA}, without differentiability assumptions in a class of two-variable means, which contains the class of \emph{Matkowski means}.

\end{abstract}

\maketitle

\section{Introduction}
We are going to use the usual notations $\N$, $\Q$, $\R$, and $\mathbb{C}$ for the sets of positive integers, rational numbers, real numbers, and complex numbers, respectively. The set of positive real numbers will be denoted by $\R_+$, that is, $\R_+:=\{x\in\R\mid x>0\}$. The identity function and the function which is identically $1$ will be denoted by the symbols $\id$ and $\mathbf{1}$, respectively.

Throughout this paper, the subset $I\subseteq\R$ will stand for a nonempty open subinterval. For a given function $h:I\to I$ and for $n\in\N\cup\{0\}$, the $n$th iterate of $h$ will be denoted by $h^{[n]}:I\to I$, where $h^{[0]}:=\id$ on $I$, furthermore $h^{[n]}(x):=h\big(h^{[n-1]}(x)\big)$ for all $x\in I$, whenever $n\in\N$.

A two-place function $M:I\times I\to\R$ will be called \emph{a two-variable mean on $I$} or, shortly, \emph{a mean on $I$} if
\Eq{M}{
\min(x,y)\leq M(x,y)\leq\max(x,y),\qquad(x,y\in I).
}
If both of the above inequalities are strict whenever $x\neq y$, then $M$ is said to be \emph{a strict mean}. We note that, by their definition, two-variable means are \emph{reflexive}, that is, we have $M(x,x)=x$ for all $x\in I$.

We say that the mean $M$ is \emph{strictly monotone} if, for all fixed $x_0,y_0\in I$, the functions
\Eq{*}{
x\mapsto M(x,y_0),\quad(x\in I)\qquad\mbox{and}\qquad
y\mapsto M(x_0,y),\quad(y\in I)
}
are strictly increasing on $I$. Observe that strictly monotone means are also strict. Finally, the mean $M$ is called \emph{symmetric} if $M(x,y)=M(y,x)$ holds for all $x,y\in I$.

The mean $M:I\times I\to\R$ is said to be a \emph{quasi-arithmetic mean} if there exists a continuous, strictly monotone function $\varphi:I\to\R$ such that
\Eq{qam}{
M(x,y)=\varphi^{-1}\bigg(\frac{\varphi(x)+\varphi(y)}2\bigg)=:\mathscr{A}_\varphi(x,y),\qquad(x,y\in I).
}
The function $\varphi$ is called \emph{the generator of the mean}. The quasi-arithmetic mean generated by $\id$ is called a two-variable arithmetic mean, which, for brevity, will be simply denoted by $\A$. In view of its definition, quasi-arithmetic means are continuous, strictly monotone, and symmetric as well.

Now, we recall some generalizations of quasi-arithmetic means, which will be crucial in our further investigations.

\emph{I. The class of Cauchy means.} A two-variable mean $M:I\times I\to\R$ is called a \emph{Cauchy mean} if there exist differentiable functions $f,g:I\to\R$ such that $0\notin g'(I)$, the ratio $f'/g'$ is invertible, and for all $x,y\in I$ with $x\neq y$, we have
\Eq{CM}{
M(x,y)=\bigg(\frac{f'}{g'}\bigg)^{-1}\bigg(\frac{f(x)-f(y)}{g(x)-g(y)}\bigg).
}
In this case $M$ is denoted by $\mathscr{C}_{f,g}$, where the functions $f$ and $g$ are called \emph{the generators of the Cauchy mean}. We note that, in view of the Cauchy Mean Value Theorem, under the above conditions, the formula in \eq{CM} indeed defines a mean on $I$, which turns out to be strict and symmetric.

To see that (at least regular) quasi-arithmetic means are contained in this class, let $\varphi:I\to\R$ be a differentiable function with non-vanishing first derivative, define further $f:=\varphi^2$ and $g:=\varphi$. Then, due to the Rolle Mean Value Theorem, the ratio $\frac{f'}{g'}=\frac{2\varphi\varphi'}{\varphi'}=2\varphi$ is invertible and, for all $x,y\in I$ with $x\neq y$, we have that
\Eq{*}{
\mathscr{C}_{f,g}(x,y)
=\big(2\varphi\big)^{-1}\bigg(\frac{\varphi^2(x)-\varphi^2(y)}{\varphi(x)-\varphi(y)}\bigg)
=\varphi^{-1}\bigg(\frac{\varphi(x)+\varphi(y)}2\bigg)=\mathscr{A}_\varphi(x,y).}

\emph{II. The class of Bajraktarević means.} This family of mean values was first investigated by the Bosnian mathematician Mahmut Bajraktarević in \cite{Baj58}. We say that a function $M:I\times I\to\R$ is a \emph{Bajraktarević mean} if there exist continuous functions $f,g:I\to\R$ such that $0\notin g(I)$, the ratio $f/g$ is invertible, and
\Eq{BM}{
M(x,y)=\bigg(\frac{f}{g}\bigg)^{-1}\bigg(\frac{f(x)+f(y)}{g(x)+g(y)}\bigg)=:\mathscr{B}_{f,g}(x,y),\qquad(x,y\in I).
}
The functions $f$ and $g$ are called \emph{the generators of the Bajraktarević mean}. Again, the definition yields that Bajraktarević means are continuous, strict, and symmetric. Furthermore, if $f:=\varphi$ is a continuous, strictly monotone function and $g=\mathbf{1}$ on $I$, then $\mathscr{B}_{f,g}=\mathscr{A}_\varphi$ on the domain $I\times I$. 

\emph{III. The class of Matkowski means.} This type of means was introduced in 2010 by the Polish mathematician Janusz Matkowski in his paper \cite{Mat10b}. A two-variable function $M:I\times I\to\R$ is called a \emph{two-variable generalized quasi-arithmetic mean} or, shortly, a \emph{Matkowski mean} if one can find continuous functions $f,g:I\to\R$ strictly monotone in the same sense such that
\Eq{MT}{
M(x,y):=(f+g)^{-1}\big(f(x)+g(y)\big),\qquad(x,y\in I).
}
The Matkowski mean generated by the pair of functions $(f,g):I\to\R^2$ will be denoted by $\mathscr{M}_{f,g}$.

It is easy to see that the expression in \eq{MT} indeed defines a mean on $I$, which is continuous and strictly monotone. If $f:=g:=\varphi$ for some continuous, strictly monotone function $\varphi:I\to\R$, then $\M_{f,g}=\A_\varphi$ on $I\times I$. Thus, under this setting, we obtain symmetric Matkowski means, however the definition \eq{MT} shows that this is not the case in general.

The following result, which characterizes quasi-arithmetic means among Matkowski means,  can be found in \cite{Mat10b}.

\Thm{MAT}{\textup{(J. Matkowski, 2010.)} Let $f,g:I\to\R$ be continuous functions, which are strictly monotone in the same sense. Then the following statements are equivalent.
\begin{enumerate}[(i)]\itemsep=1mm
\item The mean $\M_{f,g}$ is quasi-arithmetic.
\item The mean $\M_{f,g}$ is symmetric.
\item There exists a constant $c\in\R$ such that $f=g+c$ holds on $I$.
\end{enumerate}}

Now we introduce the main notion of this paper.

\section{Balancing property of means}

We say that a two-variable mean $M:I\times I\to\R$ \emph{possesses the balancing property} or, shortly, it is \emph{balanced} if, for all $x,y\in I$, we have
\Eq{B}{
M\big(M(x,u),M(u,y)\big)=u\quad\text{with}\quad u=M(x,y).}

\Rem{pe}{The prototypical example of balanced means are \emph{arithmetic means}. Indeed, if $x,y\in I$ are arbitrary and we define $u:=\A(x,y)=\frac12(x+y)$, then
\Eq{*}{
\A\big(\A(x,u),\A(u,y)\big)=\frac12\Big(\frac{3x+y}4+\frac{x+3y}4\Big)=\frac{x+y}2=\A(x,y)=u.}}

A similar calculation yields that quasi-arithmetic means are balanced as well. Instead of performing it, we draw attention to a more general phenomenon, more precisely, for later reference, we formulate the next statement about well-known inheritance properties of balancedness.

\Prp{H}{Let $M,N:I\times I\to\R$ be means having the balancing property. Then the following statements hold.
\begin{enumerate}[(i)]\itemsep=1mm
\item\label{cj} Conjugation of means preserves the balancing property: for any open subinterval $J\subseteq\R$ and continuous, strictly monotone function $\varphi:J\to I$, the mean $M_\varphi:J\times J\to\R$ defined by
\Eq{H1}{
M_\varphi(x,y):=\varphi^{-1}\big(M(\varphi(x),\varphi(y)\big)
}
is balanced.
\item\label{ft} Fitting of means preserves the balancing property: the mean $F_{M,N}:I\times I\to\R$ defined by
\Eq{H2}{
F_{M,N}(x,y):=
\begin{cases}
M(x,y)&\text{if }x\leq y,\\[1mm]
N(x,y)&\text{if }x>y
\end{cases}
}
is balanced.
\end{enumerate}
}
\begin{proof}
To prove the statement \emph{(i)}, let $x,y\in I$ be arbitrarily fixed and define $u:=M_\varphi(x,y)$. By the definition of $M_\varphi$, we obviously have
$\varphi\big(M_{\varphi}(x,u)\big)=M(\varphi(x),\varphi(u))$ and $\varphi\big(M_{\varphi}(u,y)\big)=M(\varphi(u),\varphi(y))$. Using this, then the definition of $u$, and, finally, the balancing property of $M$, we obtain that
\Eq{*}{
M_\varphi\big(M_\varphi(x,u),M_\varphi(u,y)\big)
=\varphi^{-1}\big(M\big(M(\varphi(x),\varphi(u)),M(\varphi(u),\varphi(y))\big)\big)
=\varphi^{-1}(\varphi(u))=u.
}

To show the validity of \emph{(ii)}, again, let $x,y\in I$ be any points, where, without loss of generality, we may assume that $x\neq y$. Define further $u:=F_{M,N}(x,y)$.
	
If $x<y$, then $u=M(x,y)$, furthermore $x\leq u\leq y$, where at least one of the inequalities is strict. Therefore, we have $F_{M,N}(x,u)=M(x,u)$ and $F_{M,N}(u,y)=M(u,y)$. Using these equalities and then the inequality $M(x,u)\leq M(u,y)$, we obtain that
\Eq{*}{
F_{M,N}\big(F_{M,N}(x,u),F_{M,N}(u,y)\big)
=F_{M,N}\big(M(x,u),M(u,y)\big)
=M\big(M(x,u),M(u,y)\big)=u,
}
where the very last step is due to the definition of $u$ and the balancing property of $M$.

If $y<x$, then $u=N(x,y)$ and $y\leq u\leq x$, where, again, at least one of the appearing inequalities must be strict. If $u=y$, then $F_{M,N}(x,u)=F_{M,N}(u,y)=u$, which means that the desired equality is trivially satisfied. The case $u=x$ yields the same conclusion, consequently, we may assume that $y<u<x$.

Then $\xi:=F_{M,N}(x,u)=N(x,u)$ and $\eta:=F_{M,N}(u,y)=N(u,y)$. The mean value property of $N$ implies that $\eta\leq u\leq \xi$. Furthermore, due to the definition of $u$ and the balancing property of $N$, we have $N(\xi,\eta)=u$. Therefore
\Eq{*}{
F_{M,N}&\big(F_{M,N}(x,u),F_{M,N}(u,y)\big)
=F_{M,N}\big(\xi,\eta\big)=
\begin{cases}
M\big(\xi,\eta\big)=N\big(\xi,\eta\big)=u   &\text{if }\xi=\eta,\\[1mm]
N\big(\xi,\eta\big)=u &\text{if }\xi>\eta,
\end{cases}
}
which finishes the proof.
\end{proof}

In view of \rem{pe} and statement \emph{(\ref{cj})} of \prp{H}, the balancing property of quasi-arithmetic means follows easily. To better understand which properties balancedness depends on, one can use the well-known characterization theorem of János Aczél \cite{Acz48a} as well. It states that a two-variable function $M:I\times I\to\R$ is a quasi-arithmetic mean if and only if it is strictly monotone, continuous, reflexive, symmetric, and bisymmetric, that is, for all $x,y,u,v\in I$, we have
\Eq{bys}{
M\big(M(x,y),M(u,v)\big)=M\big(M(x,u),M(y,v)\big).
}

Now, assume that $M$ is quasi-arithmetic. Applying \eq{bys} for $x,y\in I$ and for $u:=v:=M(x,y)$, then using the reflexivity and the symmetry of $M$, one can deduce \eq{B}.

Now, to demonstrate the variety of the solutions of equation \eq{B}, we give a list of a couple of some balanced means having different regularity properties.

\begin{enumerate}[(a)]
\item \emph{Solution, which is discontinuous at any point of its domain.} Let $M(x,y):=y$ whenever we have $x,y\in\Q$ and define $M(x,y):=x$ otherwise.
\item \emph{Discontinuous solution, which is homogeneous.} Let $A,B\subseteq\R_+$ be nonempty sets such that $A\cap B=\emptyset$ and $A\cup B=\R_+$. Define $M:\R_+\times\R_+\to\R$ by $M(x,y):=y$ or $M(x,y):=x$, whenever $x/y\in A$ or $x/y\in B$, respectively.
\item \emph{Continuous solutions, which are neither symmetric nor strict.} Coordinate means, that is, $M_1$ and $M_2$, defined by $M_1(x,y):=x$ and $M_2(x,y):=y$ if $(x,y)\in I\times I$.
\item \emph{Continuous, symmetric solutions, which are not strict.} Extremal means, that is, $\min$ and $\max$.
\item \emph{Continuous, strict solutions, which are not symmetric.} Let $\varphi:I\to\R$ and $\psi:I\to\R$ be continuous, strictly monotone functions such that $\{\varphi,\psi,\mathbf{1}\}$ is linearly independent over $\R$, and define $M$ by $F_{\A_\varphi,\A_\psi}$ on $I\times I$.
\item \emph{Continuous, symmetric, strict solutions.} Quasi-arithmetic means.
\end{enumerate}

A continuous, symmetric, strict solution of \eq{B}, which fails to be a quasi-arithmetic mean, is not known. We note that, to construct such an example, it is enough to focus on the continuity and strictness of the mean in question. Indeed, if $M$ is continuous, strict, and balanced, then, in view of \emph{(\ref{ft})} of \prp{H}, the symmetric mean $M^*$ defined by $M^*(x,y):=M\big(\min(x,y),\max(x,y)\big)$ inherits these properties.

It is an interesting question, beside the balancing property, what we need to assume to conclude that it is quasi-arithmetic. The first remarkable investigations in this direction are due to Georg Aumann \cite{Aum35,Aum37}. In 1935, considering the problem on the complex plain, in the paper \cite{Aum35}, it was shown that among analytic means (i.e., $M$ is reflexive, symmetric, and holomorphic on a neighborhood $G$ of a regular point $z_0\in\mathbb{C}\times\mathbb{C}$), the only solutions of equation \eq{B} are analytic quasi-arithmetic means.

Then, in 1937, turning to the real case, the author proved in \cite{Aum37} that the balancing property characterizes regular quasi-arithmetic means among Cauchy means.

In 2015, Lucio R. Berrone \cite{Ber15} introduced and investigated a generalized version of \eq{B}. Having various combinations of the conditions of symmetry, strictness and continuous differentiability, he obtained that the general solutions are quasi-arithmetic means.

Our main result will be analogous to the mentioned theorem of Aumann concerning Cauchy means. In our investigations we are going to avoid differentiability and symmetry. (In view of \thm{MAT}, the latter assumption would not be so practical.)  

\section{Auxiliary results}
In this section, let $L:I\times I\to\R$ and $R:I\times I\to\R$ be continuous, strictly monotone means. For a given $v\in I$, the continuous, strictly increasing functions $t\mapsto L(t,v)$ and $t\mapsto R(t,v)$ defined on $I$ will be denoted by $L_v$ and $R_v$, respectively. Then, particularly, $R_v$ is invertible, hence the function $\psi_v:=L_v\circ R_v^{-1}:R_v(I)\to\R$ is well-defined, as well as is continuous and strictly increasing.

For brevity, for a given point $v\in I$, let us denote the open subinterval $R_v(I)\subseteq I$ by $J_v$. Observe, that, for $v\in I$, the set $J_v$ always contains $v$.

\Lem{1}{Let $v\in I$ be arbitrarily fixed and assume that
\Eq{RL}{
R_v<L_v\quad\text{on }]-\infty,v[\,\cap\,I\qquad\text{and}\qquad
R_v>L_v\quad\text{on }I\,\cap\,]v,+\infty[\,
}
hold. Then $\psi_v(J_v)\subseteq J_v$ and, for all $\xi\in J_v\setminus\{v\}$, the sequence $\big(\psi^n_v(\xi)\big)$ converges in a strictly monotone way to the point $\psi_v(v)=v$ as $n\to\infty$.}

\begin{proof} Let $t\in J_v$ be any point with $t<v$. Then, by the definition of the interval $J_v$ and the fact that $R_v$ is injective, there uniquely exists $s\in I$ such that $R_v(s)=t$. Here, since $R_v$ is strictly increasing, we must have $s<v$. Using the first inequality in \eq{RL}, we obtain that $t=R_v(s)<L_v(s)=\psi_v(t)$, consequently, $\id<\psi_v$ holds on the interval $J_v^-:=\,]-\infty,v[\,\cap\,J_v$. A similar argument shows that $\psi_v<\id$ on $J_v^+:=J_v\,\cap\,]v,+\infty[\,$. Using the continuity of the function $\psi_v$, the fixed point property $\psi_v(v)=v$ follows. Hence,
\Eq{*}{
\psi_v(J_v)=\psi_v(J_v^-\cup\{v\}\cup J_v^+)
=\psi_v(J_v^-)\cup\psi_v(\{v\})\cup\psi_v(J_v^+)
\subseteq J_v^-\cup\{v\}\cup J_v^+=J_v.}
	
Let, finally, $\xi\in J_v\setminus\{v\}$ be arbitrarily fixed. Without loss of generality, we may assume that $\xi_0:=\xi<v$. Then, in view of the previous part of the proof, $\xi_0<\xi_1:=\psi_v(\xi_0)<\psi_v(v)=v$ follows. Applying $\psi_v$ for the point $\xi_1$ instead of $\xi_0$ and using the strict monotonicity of $\psi_v$, we get that $\xi_0<\xi_1<\xi_2:=\psi_v^{[2]}(\xi_0)<v$. By induction on the iterative power, one can obtain that the sequence $\big(\psi_v^{[n]}(\xi)\big)$ is strictly increasing and is contained in the open interval $]\xi,v[\,$. Denoting its limit by $A\in\,]\xi, v]$ and using the continuity of $\psi_v$, we have that
\Eq{*}{
\psi_v(A)
=\psi\big(\lim\limits_{n\to\infty}\psi_v^{[n]}(\xi)\big)=\lim\limits_{n\to\infty}\psi_v^{[n+1]}(\xi)=A.
}
Since $\psi_v$ is injective, $A=v$ follows. The case $v<\xi$ can be treated similarly.
\end{proof}

For a given $u\in I$, the domain of the function $v\mapsto\psi_v(u)$ will be denoted by
\Eq{dom}{D(u):=\{v\in I\mid u\in J_v=R_v(I)\}.}

\Lem{1.5}{
For any $u\in I$, the set $D(u)$ is a subinterval of $I$ containing $u$ in its interior.
}

\begin{proof}
Obviously, $u\in D(u)$, therefore $D(u)$ is nonempty. \emph{In the first step, we show that $D(u)$ cannot be a singleton}, more precisely, that $\inf D(u)<u<\sup D(u)$. Indirectly, assume that $\sup D(u)=u$ and let $r>0$ such that $u-r\in I$. Let further $u<v$ in $I$ be arbitrary. Then, by our indirect assumption and the inclusion $v\in J_v$, we must have $u\leq \inf J_v$. Using that $R_v$ is strictly increasing, we obtain that
\Eq{*}{
u\leq\inf J_v=\inf R_v(I)<R_v(u-r)=R(u-r,v).
}
Taking the limit $v\to u^+$ and using the continuity of $R$ in its second variable, we get that $u\leq R(u-r,u)<u$, because $R$ is strict. This contradiction shows that we must have $u<\sup D(u)$. The inequality $\inf D(u)<u$ can be proved similarly. 
	
\emph{In the rest of the proof we show that $D(u)$ is an interval.} Let $u'\in D(u)$ be arbitrarily fixed with $u'<u$ and let $u'<\eta<u$ be any further point of $I$. By the choice of $u'$, there exists $x\in I$, necessarily with $u<x$, such that $R_{u'}(x)=u$. Obviously, $R_{\eta}(\eta)=\eta<u$ and, using that $R$ is strictly increasing in its second variable, we have
\Eq{*}{
u=R_{u'}(x)=R(x,u')<R(x,\eta)=R_\eta(x).
}
By the Darboux Property of the function $R_{\eta}:I\to I$, we get that there exists $\xi\in\,]\eta,x[$ such that $R_{\eta}(\xi)=u$, that is, $\eta\in D(u)$. The point $\eta$ was an arbitrary element of $]u',u[\,$, consequently, we have $[u',u]\subseteq D(u)$. 

A similar argument shows that $[u,u']\subseteq D(u)$ for all $u'\in D(u)$ with $u<u'$.
\end{proof}

\Lem{2}{For any $u\in I$, the function $v\mapsto\psi_v(u)$ is continuous on $D(u)$.}

\begin{proof}
Let $u\in I$ be arbitrarily fixed and $v_0\in D(u)$. To prove that $v\mapsto\psi_v(u)=L_v\big(R_v^{-1}(u)\big)$ is continuous at $v_0$, it is enough to show that the function $v\mapsto R_v^{-1}(u)$ is continuous at $v_0$.
	
Let $(v_n)\subseteq D(u)$, different from a constant sequence, such that $v_n\to v_0$ as $n\to\infty$, furthermore define $t_n:=R^{-1}_{v_n}(u)$ if $n\in\N$ and $t_0:=R^{-1}_{v_0}(u)$. Clearly, we have to show, that $t_n\to t_0$ as $n\to\infty$.

Denote the lower limit and the upper limit of $(t_n)$ by $\alpha$ and $\beta$, respectively. Firstly, indirectly, let us assume that $\beta=\sup I$. Then, there exists a subsequence $(t_{n_k})$ of $(t_n)$ such that $t_{n_k}\to \sup I$ as $k\to\infty$. Therefore, for a given $a\in I$, there exists $k_a\in\N$ such that $a<t_{n_k}<\sup I$ whenever $k\geq k_a$. Taking the corresponding subsequence $(v_{n_k})$ of $(v_n)$, we obtain that
\Eq{*}{
R_{v_{n_k}}(a)=R(a,v_{n_k})<R(t_{n_k},v_{n_k})=u,\qquad(k\geq k_a).
}
Using the continuity of $R$ in its second variable, we get that $R(a,v_0)\leq u$ holds, where $a\in I$ was arbitrary. This leads to a contradiction whenever $a$ is chosen to be greater than $t_0$, consequently, we must have $\beta<\sup I$. A similar argument shows that $\inf I<\alpha$.
	
Now we can show the continuity of $v\mapsto R_v^{-1}(u)$ at the point $v_0$.
By the definition of $(t_n)$, we have $u=R_{v_n}(t_n)=R(t_n,v_n)$ for all $n\in\N$. Using that $R$ is continuous and strictly increasing in its variables, we get that
\Eq{*}{
u=\lim\limits_{n\to\infty}R(t_n,v_n)=\liminf\limits_{n\to\infty}R(t_n,v_n)=R\big(\liminf\limits_{n\to\infty}t_n,\liminf\limits_{n\to\infty}v_n\big)=R(\alpha,v_0)=R_{v_0}(\alpha).}
The function $R_{v_0}:I\to I$ is injective, which yields that $\alpha=t_0$. Using a similar argument, it can be shown that $\beta=t_0$ holds as well. Consequently, $(t_n)$ is convergent and tends to $t_0$ as $n\to\infty$, which finishes the proof.
\end{proof}

\section{Iteratively quasi-arithmetic means}

The solution (e) of equation \eq{B}, listed in the second section, suggests that Matkowski means, besides continuity and strict monotonicity, must have an additional property, which, together with the balancing property, yields its symmetry. Motivated by this, we introduce the following general class of means.

Let a two-variable mean $M:I\times I\to\R$ be a member of the class $\mathcal{M}(I)$ if and only if it is continuous, strictly monotone, and there exists a continuous, strictly monotone function
$\Phi:I\to\R$ such that
\Eq{IQA}{
M(x,y)=\mathscr{A}_\Phi\big(M(x,M(x,y)),M(M(x,y),y)\big),\qquad(x,y\in I).
}
If \eq{IQA} is satisfied, we are going to say that $M$ is \emph{iteratively quasi-arithmetic with respect to the generator $\Phi$}.

Roughly speaking, property \eq{IQA} says that $M$ can be obtained as a quasi-arithmetic mean of its right and left iterate. First, to demonstrate that the definition of $\mathcal{M}(I)$ is not so restrictive or artificial, we formulate and prove the following proposition.

\Prp{MIQ}{
Let $f,g:I\to\R$ be continuous functions, which are strictly monotone in the same sense. Then the Matkowski mean $\M_{f,g}$ is iteratively quasi-arithmetic with respect to the generator $\Phi:=f+g$.}

\begin{proof}
Let $x,y\in I$ be arbitrarily chosen and, for brevity, introduce the notation $u:=\M_{f,g}(x,\M_{f,g}(x,y))$ and $v:=\M_{f,g}(\M_{f,g}(x,y),y)$. Then, obviously,
\Eq{*}{
(f+g)(u)=f(x)+g\big(\M_{f,g}(x,y)\big)\quad\mbox{and}\quad
(f+g)(v)=f\big(\M_{f,g}(x,y)\big)+g(y).
}
Adding up these equalities side by side, we obtain that
\Eq{*}{
(f+g)(u)+(f+g)(v)=f(x)+(f+g)\big(\M_{f,g}(x,y)\big)+g(y)=2f(x)+2g(y).
}

Therefore we have
\Eq{*}{
\mathscr{A}_{f+g}(u,v)&=(f+g)^{-1}\bigg(\frac{(f+g)(u)+(f+g)(v)}2\bigg)=(f+g)^{-1}\bigg(\frac{2f(x)+2g(y)}2\bigg)=\M_{f,g}(x,y),
	}
	which finishes the proof.
\end{proof}

The next example shows that the class of Matkowski means is strictly contained in the class $\mathcal{M}(I)$.

\Exa{ex2}{
Let $\Phi:I\to\R$ be a continuous, strictly monotone function and $t\in\,]0,1[\,\setminus\{\frac12\}$ be arbitrarily fixed. Then the mean $K:I\to\R$ defined by
\Eq{*}{
K(x,y):=\A_\Phi^t\big(\min(x,y),\max(x,y)\big)=
\begin{cases}
\A_\Phi^t(x,y)&\text{if }x\leq y,\\[1mm]
\A_\Phi^{1-t}(x,y)&\text{if }x>y
\end{cases}
}
belongs to $\mathcal{M}(I)$, but it is not a Matkowski mean.
}
\begin{proof}
The continuity and strict monotonicity of $K$ follows easily from its definition. Now, we show that $K$ is iteratively quasi-arithmetic with respect to the function $\Phi$.

Let $x,y\in I$ be arbitrarily fixed and $u:=K(x,y)$. To avoid the trivial case, we may assume that $x\neq y$. We perform the calculation only for $x<y$, because the complementary case can be treated similarly. Then $u=\A_\Phi^t(x,y)$ and, obviously, $x<u<y$, therefore
\Eq{*}{
\A_\Phi\big(K(x,u),K(u,y)\big)
&=\A_\Phi\big(\A_\Phi^t(x,u),\A_\Phi^t(u,y)\big)
=\Phi^{-1}\bigg(\frac{t\Phi(x)+\Phi(u)+(1-t)\Phi(y)}{2}\bigg)\\
&=\Phi^{-1}\big(t\Phi(x)+(1-t)\Phi(y)\big)
=\A_\Phi^t(x,y)=K(x,y).
}

Finally, if $K$ were a Matkowski mean, then, by its symmetry, in view of \thm{MAT}, it would be a quasi-arithmetic mean. On the other hand it can be shown that $K$ is not bisymmetric. Indeed, applying equation \eq{bys} for any points $x<u<y$ in $I$ and $v:=y$ under $M=K$, we obtain
\Eq{*}{
(1-t)(2t-1)\big(\Phi(y)-\Phi(u)\big)=0\qquad\text{or}\qquad
t(2t-1)\big(\Phi(u)-\Phi(x)\big)=0,
}
provided that $K(x,y)\leq u$ or $K(x,y)>u$, respectively. In view of the definition of $t$ and the injectivity of $\Phi$, both are impossible.
\end{proof}

Now we are going to formulate and prove an extension theorem, which will be crucial proving our main result. We note that similar investigations can be found in the paper \cite{Mak05} of Gyula Maksa related to more general functions, namely quasi-sums. In \cite{Mak05}, the author proves that if a function is a quasi-sum on some rectangular neighborhood of any point of its domain, that is, it is a local quasi-sum, then it can be written as a quasi-sum on its entire domain. Roughly speaking, the property is that the quasi-sum is localizable.

In our situation, we assume that our two-place function is a quasi-arithmetic mean on some neighborhood of any point of the diagonal of its domain, that is, diagonally locally quasi-arithmetic, and that it is iteratively quasi-arithmetic on its entire domain. For our purposes, it is enough to assume that the generator functions of the quasi-arithmetic means mentioned here are the same.

\Thm{3}{Let $M\in\mathcal{M}(I)$ and $\Phi:I\to\R$ be a continuous, strictly monotone function, for which \eq{IQA} holds with $M$. If, for all $p\in I$, there exists an open neighborhood $U_p\subseteq I$ of $p$ such that $M=\A_\Phi$ on the rectangle $U_p\times U_p$, then $M=\A_\Phi$ on $I\times I$.
}

\begin{proof}
Let $p\in I$ be arbitrarily fixed and $U_p\subseteq I$ the corresponding neighborhood of $p$. If $U_p=I$ then we are done. \emph{Hence, indirectly, let us assume that the intersection $\{a_0:=\inf U_p, b_0:=\sup U_p\}\cap I$ is nonempty, say we have $a_0\in I$.} The case $b_0\in I$ can be treated similarly. We may also assume that $U_p$ is maximal in $I$, that is, for any open subinterval $U_p\subsetneq J\subseteq I$, there exist $(x,y)\in J\times J$ such that $M(x,y)\neq\mathscr{A}_\Phi(x,y)$.

By the continuity of $M$ in its variables, it follows that $a_0\in U_p$. Then, in view of our condition, there exists a neighborhood $U_{a_0}$ of $a_0$ such that $M=\mathscr{A}_\Phi$ on the product $U_{a_0}\times U_{a_0}$. Again, we may assume that $U_{a_0}$ is maximal in $I$. Observe further that $b_1:=\sup U_{a_0}<b_0$ must hold. Otherwise, the interval $U_p$ were expandable, contradicting its maximality. Therefore $b_1\in U_{a_0}\cap I$ holds as well.

Due to the definition of the interval $U_p$ and our indirect assumption $a_0\in I$, for all $\varepsilon>0$, there exists a point
$x\in A_\varepsilon:=\,]a_0-\varepsilon,a_0[\,\cap\,U_{a_0}$ such that the set
\Eq{*}{
Y(x):=\{y\in U_p\mid (M-\mathscr{A}_\Phi)(x,y)\neq 0\}
}
is nonempty. Moreover, by the continuity of $M$ in its second variable, $Y(x)$ is open in $U_p$. The definition of $x$ also implies that we must have $Y(x)\subseteq\,]b_1,+\infty[\,\cap\,U_p$.

Now, temporarily, let $\varepsilon>0$ and $x\in A_\varepsilon$ be arbitrarily fixed so that $Y(x)$ is nonempty. We are going to show that, for all $y\in Y(x)$, the inclusion $M(x,y)\in I\,\setminus\,]a_0,b_1[\,$ holds.

Indirectly, assume that this is not the case, that is, there exists $y\in Y(x)$, such that $a_0<M(x,y)<b_1$. Then, particularly, $M(x,y)\in U_{a_0}\cap U_p$, therefore we obtain that
\Eq{*}{
M\big(x,M(x,y)\big)=\mathscr{A}_\Phi\big(x,M(x,y)\big)
\qquad\mbox{and}\qquad
M\big(M(x,y),y\big)=\mathscr{A}_\Phi\big(M(x,y),y\big).
}
Using that $M$ is iteratively quasi-arithmetic with respect to $\Phi$, we get that
\Eq{*}{
M(x,y)
=\mathscr{A}_\Phi\big(M(x,M(x,y)),M(M(x,y),y)\big)=
\Phi^{-1}\bigg(\frac{\Phi(x)+2\Phi(M(x,y))+\Phi(y)}{4}\bigg).
}
Applying $\Phi$ on both sides, then expressing $M(x,y)$, it follows that $M(x,y)=\mathscr{A}_\Phi(x,y)$, which contradicts the definition of $y$.

Motivated by this, define
\Eq{*}{
Y^-(x):=\{y\in Y(x)\mid M(x,y)\leq a_0\}
\qquad\text{and}\qquad
Y^+(x):=\{y\in Y(x)\mid b_1\leq M(x,y)\}.
}
Then, in view of the previous argument, $Y(x)=Y^-(x)\cup Y^+(x)$. \emph{To get a contradiction, in the rest of the proof we show that both of the sets $Y^-(x)$ and $Y^+(x)$ must be empty.}

Let $(x_n)\subseteq I$ be a sequence with $x_n\in A_{1/n}$ and $Y(x_n)\neq\emptyset$ for all $n\in\N$. Then, for a given $n\in\N$, let $y_n\in Y(x_n)$ be arbitrarily fixed. In view of the previous argument, at least one of the inclusions $y_n\in Y^-(x_n)$ and $y_n\in Y^+(x_n)$ holds for infinitely many indices. Hence, without loss of generality, we may assume that $y_n\in Y^-(x_n)$ holds for all $n\in\N$.

Define $u_n:=M(x_n,y_n)$ whenever $n\in\N$. By the definition of $(x_n)$, we get that $(u_n)$ is convergent having the limit $a_0$. Let $y^*:=\limsup_{n\to\infty} y_n$. Obviously, we have $b_1\leq y^*$. We claim that $y^*<b^*:=\sup I$ also holds. Otherwise, there were a subsequence $(y_{n_k})$ of $(y_n)$ such that $y_{n_k}\to b^*$ as $k\to\infty$. If $a\in\,]a_0,b^*[\,$ is arbitrarily fixed and $k_0\in\N$ such that $a<y_{n_k}<b^*$ whenever $k\geq k_0$, then, using that $M$ is strictly monotone in its second variable and then the inclusion $(y_n)\subseteq \bigcup_{n\in\N}Y^-(x_n)$, we have that
\Eq{*}{
M(x_{n_k},a)<M(x_{n_k},y_{n_k})\leq a_0,\qquad(k\geq k_0).
}
Due to the continuity of $M$ in its first variable, we obtain finally that $M(a_0,a)\leq a_0$, contradicting that $M$ is a strict mean. Hence $y^*\in I$, thus we can substitute it into $M$ and we get that
\Eq{*}{
a_0
=\lim\limits_{n\to\infty}M(x_n,y_n)
=\limsup\limits_{n\to\infty}M(x_n,y_n)
=M\big(\limsup\nolimits_{n\to\infty} x_n, \limsup\nolimits_{n\to\infty} y_n\big)
=M(a_0,y^*).
}
Consequently, $y^*=a_0$ must hold, which, since $a_0<b_1\leq y^*$, is a contradiction again.

Using a similar argument, one can show that $(y_n)\subseteq \bigcup_{n\in\N}Y^+(x_n)$ is impossible as well. This contradiction was caused by the assumption $a_0\in I$, consequently we must have $a_0=\inf I$. As we mentioned before, the equality $b_0=\sup I$ can be proved similarly.
\end{proof}

\section{The main result}

Now we can formulate and prove our main results.

\Thm{2}{
If $M\in\mathcal{M}(I)$ is balanced and $\Phi:I\to\R$ is a continuous, strictly monotone function for which \eq{IQA} holds with $M$, then, for all $p\in I$, there exists an open neighborhood $U_p\subseteq I$ of the point $p$ such that $M=\A_\Phi$ holds on $U_p\times U_p$.}
\begin{proof}
Define the means $L:I\times I\to\R$ and $R:I\times I\to\R$ by	
\Eq{*}{
L(s,t):=M(M(s,t),t)
\qquad\text{and}\qquad
R(s,t):=M(s,M(s,t)).
}

Then $L$ and $R$ are continuous and strictly monotone. Let $p\in I$ be arbitrarily fixed. \emph{In the first part of the proof we construct the proper open neighborhood $U_p$ of the point $p$.} In view of \lem{1.5}, the set $D(p)$, defined in \eq{dom}, is a subinterval of $I$ containing $p$ in its interior. Hence, let $v^*\in D(p)$ be arbitrary with $p<v^*$. We claim that there exists $u^*<p$ in the interval $J_{v^*}$ such that $p$ does not belong to the orbit $\{\psi_{v^*}^{[n]}(u^*)\mid n\in\N\}\subseteq J_{v^*}$, where $J_{v^*}:=R_{v^*}(I)$ and $\psi_{v^*}:=L_{v^*}\circ R_{v^*}^{-1}$.

Indirectly, assume that this is not the case, and let $u\in J_{v^*}$ with $u<p$ be any point. Then, by our indirect assumption, there exists $k\in\N$ such that $p=\psi_{v^*}^{[k]}(u)$. By \lem{1} and as $u<p<v^*$, the sequence $(\psi_{v^*}^{[n]}(u))$ is strictly increasing, hence $u\leq \psi_{v^*}^{[k-1]}(u)<p$. Now, pick a point $u'$ arbitrarily with $\psi_{v^*}^{[k-1]}(u)<u'<p$. Then, again, by our indirect assumption, there exists an index $m\in\N$ such that $\psi_{v^*}^{[m]}(u')=p$.  Therefore we have \Eq{egy}{\psi_{v^*}^{[k]}(u)=\psi_{v^*}^{[m]}(u').}
The function $\psi_{v^*}^{[n]}:J_{v^*}\to J_{v^*}$ is invertible for each fixed $n\in\N$, consequently $k$ and $m$ cannot be equal. 

If $k<m$, then, applying the inverse function of $\psi^{[k]}_{v^*}$ on both sides of the equality \eq{egy}, we get that $u=\psi_{v^*}^{[m-k]}(u')$. This means that $u$ belongs to the orbit of $u'$. Since $u<u'$ and the sequence $(\psi_{v^*}^{[n]}(u'))$ is strictly increasing, this is impossible. Thus we must have $m<k$. Similarly, applying the inverse of $\psi_{v^*}^{[m]}$ on both sides of \eq{egy}, we obtain that $u'$ belongs to the orbit of $u$. This, contradicts the definition of $u'$.

Hence, let $u^*<p$ in $J_{v^*}$ and $k\in\N$ such that $a:=\psi_{v^*}^{[k]}(u^*)<p<\psi_{v^*}^{[k+1]}(u^*)=:b$, and let us define $U_p:=\,]a,b[\,\subseteq J_{v^*}$. By \lem{1}, the inclusion $U_p\subseteq\, ]-\infty,v^*[\,\cap\,J_{v^{*}}$ holds as well. \emph{In the remaining part of the proof we show that $M$ is a quasi-arithmetic mean on the open interval $U_p$ generated by the function $\Phi$.}

Let $x,y\in U_p\subseteq J_{v^*}$ be any points. We may assume that $x<y$, because the case $y<x$ can be treated similarly. By the inclusion $x\in J_{v^*}$, it follows that $v^*\in D(x)$. In view of \lem{1.5}, the interval $[x,v^*]$ is contained in $D(x)$, hence, by \lem{2}, the function $t\mapsto\psi_t(x)$ is continuous on the interval $[x,v^*]$. Obviously, $\psi_{x}(x)=x$, furthermore, in view of the relation $a<x$, we also have $b=\psi_{v^*}(a)<\psi_{v^*}(x)$. By the Darboux property of $t\mapsto \psi_t(x)$ on the interval $[x,v^*]$, one can find $v_0\in\,]x,v^*[\,\subseteq D(x)$ such that $\psi_{v_0}(x)=y$. Having the point $v_0\in D(x)$, let $u_0\in I$ be the only element for which $R_{v_0}(u_0)=x$ holds. Thus
\Eq{*}{
y=\psi_{v_0}(x)=\big(L_{v_0}\circ R_{v_0}^{-1}\big)(x)=L_{v_0}\big(R_{v_0}^{-1}(x)\big)=L_{v_0}(u_0).
}
Using the balancing property of $M$ at the pair $(u_0,v_0)$, and then, at the same pair, that $M$ is iteratively quasi-arithmetic, we get that
\Eq{*}{M(x,y)&=
M\big(R_{v_0}(u_0),L_{v_0}(u_0)\big)
=M\big(M(u_0,M(u_0,v_0)),M(M(u_0,v_0),v_0)\big)\\
&=M(u_0,v_0)=\A_\Phi\big(M(u_0,M(u_0,v_0)),M(M(u_0,v_0),v_0)\big)=\mathscr{A}_\Phi\big(R_{v_0}(u_0),L_{v_0}(u_0)\big)=\mathscr{A}_{\Phi}(x,y).
}
The points $x$ and $y$ were arbitrary elements of $U_p$, which yields that $M=\A_\Phi$ indeed holds on $U_p\times U_p$.
\end{proof}

\Thm{Main}{
Let $M\in\mathcal{M}(I)$ and $\Phi:I\to\R$ be a continuous, strictly monotone function for which \eq{IQA} holds with $M$. The mean $M$ enjoys the balancing property if and only if it is a quasi-arithmetic mean on $I\times I$ generated by the function $\Phi$.}

\begin{proof}
If $M$ is the quasi-arithmetic mean on $I\times I$ generated by $\Phi$ then it is obviously balanced.

Assume that $M$ enjoys the balancing property. Then, by \thm{2}, for all $p\in I$, there exists a rectangular neighborhood of $(p,p)$ in $I\times I$, where $M$ can be written as $\A_\Phi$. Applying \thm{3}, we obtain that the equality $M=\A_\Phi$ holds on the entire domain $I\times I$.
\end{proof}

Finally we formulate the analogue of Aumann's result concerning Cauchy means.

\Cor{C}{
Let $f,g:I\to\R$ be continuous functions, which are strictly monotone in the same sense. The Matkowski mean $\M_{f,g}$ enjoys the balancing property if and only if it is a quasi-arithmetic mean generated by the function $f+g$.
}

\begin{proof}
The statement directly follows from \prp{MIQ} and \thm{Main}.
\end{proof}

\Rem{*}{It is an interesting question what the solutions of equation \eq{B} are in the class of Bajraktarević means. To prove that they are necessarily quasi-arithmetic, we need a different argument, since these means may not be contained in the class $\mathcal{M}(I)$.}


\begin{thebibliography}{10}
\bibitem{Acz48a}
J.~Aczél.
\newblock On mean values.
\newblock {\em Bull. Amer. Math. Soc.}, 54:392–400, 1948.

\bibitem{Aum35}
G.~Aumann.
\newblock Aufbau von mittelwerten mehrerer argumente ii., (analytische
mittelwerte.).
\newblock {\em Math. Ann.}, 111(1):713--730, 1935.

\bibitem{Aum37}
Georg Aumann.
\newblock Vollkommene {F}unktionalmittel und gewisse
{K}egelschnitteigenschaften.
\newblock {\em J. Reine Angew. Math.}, 176:49--55, 1937.

\bibitem{Baj58}
M.~Bajraktarević.
\newblock Sur une équation fonctionnelle aux valeurs moyennes.
\newblock {\em Glasnik Mat.-Fiz. Astronom. Društvo Mat. Fiz. Hrvatske Ser.
	II}, 13:243–248, 1958.

\bibitem{Ber15}
Lucio~R. Berrone.
\newblock The Aumann functional equation for general weighting procedures.
\newblock {\em Aequationes Mathematicae}, 89(4):1051--1073, Aug 2015.

\bibitem{Mak05}
Gy. Maksa.
\newblock {Quasisums and generalized associativity}.
\newblock {\em Aequationes Math.}, 69(1-2):6–27, 2005.

\bibitem{Mat10b}
J.~Matkowski.
\newblock Generalized weighted and quasi-arithmetic means.
\newblock {\em Aequationes Math.}, 79(3):203–212, 2010.
\end{thebibliography}

\end{document}